\documentclass[11.05pt,a4paper]{amsart}
\usepackage{amssymb,amsmath}
\usepackage{latexsym}
\usepackage{amsthm,amsfonts,amssymb,mathrsfs}
\usepackage{rotating}
\usepackage[leqno]{amsmath}
\usepackage{xspace}
\usepackage[all]{xy}
\usepackage{longtable}

\headsep=1cm \numberwithin{equation}{section}
\newtheorem{theorem}{Theorem}[section]
\newtheorem{lemma}[theorem]{Lemma}
\newtheorem{proposition}[theorem]{Proposition}
\newtheorem{corollary}[theorem]{Corollary}

\theoremstyle{definition}

\theoremstyle{remark}

\newtheorem{example}[theorem]{Example}

\newcommand{\im}{\operatorname{im}}

\newcommand{\cd}{\operatorname{cd}}

\newcommand{\fd}{\operatorname{fd}}
\newcommand{\Rfd}{\operatorname{Rfd}}

\newcommand{\Gfd}{\operatorname{Gfd}}

\newcommand{\Cone}{\operatorname{Cone}}
\newcommand{\Tor}{\operatorname{Tor}}
\newcommand{\Hom}{\operatorname{Hom}}

\newcommand{\width}{\operatorname{width}}

\newcommand{\vpl}{\operatornamewithlimits{\varprojlim}}

\newcommand{\lo}{\longrightarrow}
\newcommand{\fm}{\frak{m}}

\newcommand{\fa}{\frak{a}}

\newenvironment{prf}[1][Proof]{\begin{proof}[\bf #1]}{\end{proof}}

\begin{document}

\author[F. Mohammadi Aghjeh Mashhad and K. Divaani-Aazar]{Fatemeh Mohammadi
Aghjeh Mashhad and Kamran Divaani-Aazar}
\title[Local homology and Gorenstein flat modules]
{Local homology and Gorenstein flat modules}

\address{F. Mohammadi Aghjeh Mashhad, Science and Research Branch, Islamic Azad
University, Tehran, Iran.} \email{mohammadi\_fh@yahoo.com}

\address{K. Divaani-Aazar, Department of Mathematics, Az-Zahra
University, Vanak, Post Code 19834, Tehran, Iran.}
\email{kdivaani@ipm.ir}



\begin{abstract} Let $R$ be a commutative Noetherian ring, $\fa$ an
ideal of $R$ and $\mathcal{D}(R)$ denote the derived category of
$R$-modules. We investigate the theory of local homology in
conjunction with Gorenstein flat modules. Let $X$ be a homologically
bounded to the right complex and $Q$ a bounded to the right complex
of Gorenstein flat $R$-modules such that $Q$ and $X$ are isomorphic
in $\mathcal{D}(R)$. We establish a natural isomorphism ${\bf
L}\Lambda^{\fa}(X)\simeq \Lambda^{\fa}(Q)$ in $\mathcal{D}(R)$ which
immediately asserts that $\sup {\bf L}\Lambda^{\fa}(X)\leq \Gfd_RX$.
This isomorphism yields several consequences. For instance, in the
case $R$ possesses a dualizing complex, we show that  $\Gfd_R {\bf
L}\Lambda^{\fa}(X)\leq \Gfd_RX$. Also, we establish a criterion for
regularity of Gorenstein local rings.
\end{abstract}

\maketitle

\section{Introduction and Prerequisites}

Throughout this paper $R$ is a commutative Noetherian ring and
$\mathcal{D}(R)$ denotes the derived category of $R$-modules.  The
full subcategory of homologically bounded complexes is denoted by
$\mathcal{D}_{\Box}(R)$ and that of complexes homologically bounded
to the right (resp. left) is denoted by $\mathcal{D}_{\sqsupset}(R)$
(resp. $\mathcal{D}_{\sqsubset}(R)$). Also,
$\mathcal{D}_{\Box}^f(R)$ stands for the full subcategory of
homologically bounded complexes with finitely generated homology
modules. We use the symbol $\simeq$ for denoting isomorphisms in the
category $\mathcal{D}(R)$. For any complex $X$ in
$\mathcal{D}_{\sqsupset}(R)$ (resp. $\mathcal{D}_{\sqsubset}(R)$),
there is a bounded to the right (resp. left) complex $U$ of
projective (resp. injective) $R$-modules such that $U\simeq X$. A
such complex $U$ is called a projective (resp. injective) resolution
of $X$. We say that a homologically bounded complex $X$ has finite
projective (resp. injective) dimension if in $\mathcal{D}(R)$ it is
isomorphic to a bounded complex of projective (resp. injective)
$R$-modules. The left derived tensor product functor
$-\otimes_R^{{\bf L}}\sim$ is computed by taking a projective
resolution of the first argument or of the second one. Also, the
right derived homomorphism functor ${\bf R}\Hom_R(-,\sim)$ is
computed by taking a projective resolution of the first argument or
by taking an injective resolution of the second one.

Let $\fa$ be an ideal of $R$ and $\mathcal{C}_0(R)$ denote the full
subcategory of $R$-modules. It is known that the $\fa$-adic
completion functor
$$\Lambda^{\fa}(-)=\underset{n}{\vpl}(R/\fa^n\otimes_R-):
\mathcal{C}_0(R)\to \mathcal{C}_0(R)$$ is not right exact in
general. The left derived functor of $\Lambda^{\fa}(-)$ exists in
$\mathcal{D}(R)$, and so for any complex $X\in
\mathcal{D}_{\sqsupset}(R)$, the complex ${\bf L}\Lambda^{\fa}(X)\in
\mathcal{D}_{\sqsupset}(R)$ is defined by ${\bf
L}\Lambda^{\fa}(X):=\Lambda^{\fa}(P)$, where $P$ is a (every)
projective resolution of $X$. Let $X\in \mathcal{D}_{\sqsupset}(R)$.
For any integer $i$, the $i$-th local homology module of $X$ with
respect to $\fa$ is defined by $H_i^{\fa}(X):=H_i({\bf
L}\Lambda^{\fa}(X)).$ The study of local homology modules was
initiated by Matlis \cite{M}. Then it was continued by many authors,
see e.g. \cite{Si}, \cite{GM}, \cite{LLT}, \cite{Sc} and \cite{Fr}.
Let $\Check{C}(\underline{\fa})$ denote the $\Check{C}$ech complex
of $R$ on a set $\underline{\fa}$ of generators of $\fa$. By
\cite[(0.3),aff,p.4]{LLT} (see also \cite[Section 4]{Sc} for
corrections), $${\bf L}\Lambda^{\fa}(X)\simeq {\bf
R}\Hom_R(\Check{C} (\underline{\fa}),X).$$ By using this isomorphism
Frankild \cite[Theorem 2.11]{Fr} proved that $\inf {\bf
L}\Lambda^{\fa}(X)=\width_R(\fa,X)$, where $\width_R(\fa,X):=\inf
(R/\fa\otimes_R^{\bf L}X)$. Finding a good upper bound for $\sup
{\bf L}\Lambda^{\fa}(X)$ was considered in \cite{Sc} and \cite{Fr}.
The study of connections between Gorenstein injective modules and
local cohomology modules was started by Sazeedeh \cite{Sa}. Here, we
investigate connections between Gorenstein flat modules and local
homology modules. The notion of Gorenstein flat modules was
introduced by Enochs, Jenda and Torrecillas in \cite{EJT}. An
$R$-module $T$ is said to be Gorenstein flat if there exists an
exact complex $F$ of flat $R$-modules such that $T \cong \im(F_0\lo
F_{-1})$ and $F\otimes_RI$ is exact for all injective $R$-modules
$I$. The Gorenstein flat dimension of $X$ is defined by
$$\begin{array}{llll}\Gfd_RX&:=\inf\{\sup\{l\in \mathbb{Z}|Q_l\neq 0\}|
Q \  \ \text{is a bounded to the right complex of} \\
&\text{Gorenstein flat R-modules and} \   Q\simeq X \}.
\end{array}$$
For more details on the theory of Gorenstein homological dimensions
for complexes, we refer the reader to \cite{C}.

Let $T$ be a Gorenstein flat $R$-module and $X\in
\mathcal{D}_{\sqsupset}(R)$. We show that $T$ is
$\Lambda^{\fa}$-acyclic and $H^{\fa}_0(T)\cong \Lambda^{\fa}(T)$.
Using this, we prove that if $Q$ is a bounded to the right complex
of Gorenstein flat $R$-modules such that $Q\simeq X$, then ${\bf
L}\Lambda^{\fa}(X)\simeq\Lambda^{\fa}(Q)$, in particular $\sup {\bf
L}\Lambda^{\fa}(X)\leq \Gfd_RX$. We deduce several applications. We
show that the large restricted flat dimension of $\Lambda^{\fa}(T)$
is zero. Now, assume that $R$ possesses a dualizing complex. Then,
we prove that $\Lambda^{\fa}(T)$ is Gorenstein flat. Also, we
establish the inequality $\Gfd_{R}{\bf L}\Lambda^{\fa}(X)\leq
\Gfd_RX$, which improves \cite[Theorem 5.10 b)]{CFH}. Suppose that
$X$ is homologically bounded and let $Y\in \mathcal{D}_{\Box}^f(R)$
be a non-exact complex. If either projective or injective dimension
of $Y$ is finite, then we show that
$$\sup {\bf L}\Lambda^{\fa}(X\otimes_R^{{\bf L}}Y)\leq \Gfd_RX+\sup
Y.$$ Finally, we prove that a Gorenstein local ring $(R,\fm)$ is
regular if and only if the $\fm$-adic completion of any Gorenstein
flat $R$-module is flat.

\section{the results}

For proving our main result, we need the following three lemmas. For
an ideal $\fa$ of $R$, $\cd_{\fa}(R)$ denotes the supremum of $i$'s
such that $i$-th local cohomology module of $R$ with respect to
$\fa$ is nonzero.

\begin{lemma} Let $\fa$ be an ideal of the Noetherian ring $R$ and $X\in
\mathcal{D}_{\Box}(R)$. Then $$\sup {\bf L}\Lambda^{\fa}(X)\leq \sup
X+\cd_{\fa}(R).$$
\end{lemma}

\begin{prf} Let $\underline{\fa}$ be a set of generators of $\fa$. Then,
we have ${\bf L}\Lambda^{\fa}(X)\simeq {\bf
R}\Hom_R(\Check{C}(\underline{\fa}),X)$. Hence, \cite[Proposition
A.4.6]{C} implies that $$\sup {\bf L}\Lambda^{\fa}(X)\leq \sup
X-\inf \Check{C}(\underline{\fa})=\sup X+\cd_{\fa}(R).$$
\end{prf}

Let $\fa$ be an ideal of $R$. Simon \cite[Section 5.1]{Si}
investigated the class $C_{\fa}$ of $R$-modules $M$ for which
$H^{\fa}_i(M)=0$ for all $i>0$, and such that the natural
homomorphism $H^{\fa}_0(M)\lo \Lambda^{\fa}(M)$ is an isomorphism.
By \cite[Corllary 4.5]{M} every flat $R$-module belongs to
$C_{\fa}$. Next, we improves this result of Matlis by showing that
every Gorenstein flat $R$-module belongs to $C_{\fa}$.

\begin{lemma} Let $\fa$ be an ideal of the Noetherian ring $R$ and $Q$
a Gorenstein flat $R$-module.
\begin{enumerate}
\item[i)] $Q$ is $\Lambda^{\fa}$-acyclic.
\item[ii)] There is a natural $R$-isomorphism $H^{\fa}_0(Q)\cong
\Lambda^{\fa}(Q)$.
\end{enumerate}
\end{lemma}

\begin{prf} i) There exists an exact sequence
$$F= \cdots \to F_1\to F_0\to F_{-1}\to F_{-2}\to \cdots$$
of flat $R$-modules such that $Q\cong \im(F_0\to F_{-1})$. Set
$L_0:=Q$ and $L_i:=\im(F_i\to F_{i-1})$ for all $i<0$. For each
$i<0$, the exact sequence
$$0\to L_{i+1}\to F_i\to L_i\to 0,$$ yields the following long exact
sequence of local homology modules
$$\cdots \to H^{\fa}_{j+1}(F_i)\to H^{\fa}_{j+1}(L_i)\to
H^{\fa}_j(L_{i+1})\to H^{\fa}_j(F_i)\to \cdots.$$ The argument of
\cite[Corollary 4.5]{M} yields that any flat $R$-module is
$\Lambda^{\fa}$-acyclic. Hence, we conclude the isomorphisms
$H^{\fa}_{j}(L_{i+1})\cong H^{\fa}_{j+1}(L_i)$ for all $i<0$ and all
$j\geq 1$. Let $n:=\cd_{\fa}(R)$. Then by Lemma 2.1, one has
$$H^{\fa}_j(Q)\cong H^{\fa}_{j+1}(L_{-1})\cong \dots
\cong H^{\fa}_{j+n}(L_{-n})=0,$$  for all $j>0$.

ii) By \cite[Section 5.1]{Si},  there is a natural transformation of
functors $\xi : H^{\fa}_0(\cdot) \to  \Lambda^{\fa}(\cdot)$, which
is such that $\xi_N$ is surjective for all $R$-modules $N$. The
proof of \cite[Corollary 4.5]{M} implies that $\xi_F$ is an
isomorphism for any flat $R$-module $F$. From the definition of
Gorenstein flat $R$-modules, one can construct an exact sequence
$0\to Q\overset{f} \to F\overset{g}\to N\to 0$ of $R$-modules and
$R$-homomorphisms in which $F$ is flat and $N$ is Gorenstein flat.
By i), we deduce the following exact sequence  $$0\to H^{\fa} _0(Q)
\overset{H^{\fa}_0(f)}\longrightarrow H^{\fa}_0 (F)
\overset{H^{\fa}_0(g)}\longrightarrow H^{\fa}_0 (N)\to 0.$$ Now,
from the commutative square

\begin{equation*}
\begin{split}
\xymatrix {  H^{\fa}_0(Q) \ar[d]^-{\xi_Q} \ar[rr]^-{H_0^{\fa}(f)} &
& H^{\fa}_0(F)
\ar[d]^-{\xi_{F}} \\
\Lambda^{\fa}(Q) \ar[rr]^-{\Lambda^{\fa}(f)} & & \Lambda^{\fa}(F),}
\end{split}
\end{equation*}
it becomes clear that $\xi_Q$ is an isomorphism, as required.
\end{prf}

Next, we record the following immediate corollary of Lemma 2.2.

\begin{corollary} Let $\fa$ be an ideal of the Noetherian ring $R$. The
functor $\Lambda^{\fa}(-)$ is exact on the full subcategory of
Gorenstein flat $R$-modules.
\end{corollary}

The following useful lemma is well-known, and so we skip its proof.

\begin{lemma} Let $T:\mathcal{C}_0(R)\to \mathcal{C}_0(R)$ be a covariant
additive functor. Any morphism of complexes $\alpha:X\to Y$ yields
an isomorphism of complexes $\phi_{\alpha}:\Cone(T(\alpha))\lo
T(\Cone(\alpha))$.
\end{lemma}

Next, we present our main result.

\begin{theorem} Let $\fa$ be an ideal of the Noetherian ring $R$.
Let $X\in \mathcal{D}_{\sqsupset}(R)$ and $Q$ a bounded to the right
complex of Gorenstein flat $R$-modules such that $Q\simeq X$. Then
${\bf L}\Lambda^{\fa}(X)\simeq\Lambda^{\fa}(Q)$, and so
$H^{\fa}_i(X)=H_i(\Lambda^{\fa}(Q))$ for all $i\in \mathbb{Z}$. In
particular, $\sup {\bf L}\Lambda^{\fa}(X)\leq \Gfd_RX$.
\end{theorem}

\begin{prf}  Let $P$ be a projective resolution of $X$. Then
$P\simeq Q$, and hence  \cite[1.1.P and 1.4.P]{AF} and
\cite[A.4.1]{C} yield the existence of a quasi-isomorphism
$\alpha:P\to Q$. Now, $\Cone(\alpha)$ is an exact bounded to the
right complex of Gorenstein flat $R$-modules. By splitting
$\Cone(\alpha)$ into short exact sequences and using
\cite[Proposition 3.12]{H} and Corollary 2.3, we see that
$\Lambda^{\fa}(\Cone(\alpha))$ is exact, and so by Lemma 2.4,
$\Cone(\Lambda^{\fa}\alpha)$ is also exact. Therefore
$\Lambda^{\fa}(\alpha):\Lambda^{\fa}(P)\to \Lambda^{\fa}(Q)$ is a
quasi-isomorphism, and so $${\bf L}\Lambda^{\fa}(X)\simeq
\Lambda^{\fa}(P)\simeq  \Lambda^{\fa}(Q).$$
\end{prf}

Corollaries 2.8 and 2.10 are the main applications of this theorem.
To prove Corollary 2.8, we need a couple of lemmas. The first lemma
slightly improves \cite[1.10]{FI}. Recall that for a complex $X\in
\mathcal{D}_{\sqsupset}(R)$, any bounded to the right complex $F$
such that $F$ consists of flat $R$-modules and there exists a
quasi-isomorphism $\alpha:F\lo X$ is called a flat resolution of
$X$.

\begin{lemma} Let ${\fa}$ be an ideal of the Noetherian ring $R$
and $X,Y\in \mathcal{D}_{\sqsupset}(R)$. Let $Q$ be a bounded to the
right complex of Gorenstein flat R-modules such that $Q\simeq X$ and
$F$ a flat resolution of $Y$. Then ${\bf
L}\Lambda^{\fa}(X\otimes_R^{{\bf L}}Y)\simeq
\Lambda^{\fa}(Q\otimes_RF).$ Moreover, if $X$ is homologically
bounded and all homology modules of $Y$ are finitely generated, then
${\bf L}\Lambda^{\fa}(X\otimes_R^{{\bf L}}Y)\simeq {\bf
L}\Lambda^{\fa}(X)\otimes_R^{{\bf L}}Y.$
\end{lemma}

\begin{prf} By \cite[Ascent table II a)]{CH} $Q\otimes_RF$ is a complex
of Gorenstein flat $R$-modules. Hence, as the complex $Q\otimes_RF$
is bounded to the right and $Q\otimes_RF\simeq X\otimes_R^{{\bf
L}}Y$, Theorem 2.5 implies that ${\bf
L}\Lambda^{\fa}(X\otimes_R^{{\bf L}}Y)\simeq
\Lambda^{\fa}(Q\otimes_RF).$

Now, assume that all homology modules of $Y$ are finitely generated.
Since, by \cite[5.8]{CFH}, $\Check{C}(\underline{\fa})$ has finite
projective dimension, \cite[Proposition 2.2 vi)]{CH} yields that
$${\bf L}\Lambda^{\fa}(X\otimes_R^{{\bf L}}Y)\simeq {\bf
R}\Hom_R(\Check{C}(\underline{\fa}),X\otimes_R^{{\bf L}}Y)\simeq
{\bf R}\Hom_R(\Check{C}(\underline{\fa}),X)\otimes_R^{{\bf
L}}Y\simeq {\bf L}\Lambda^{\fa}(X)\otimes_R^{{\bf L}}Y.$$
\end{prf}

Recall that the large restricted flat dimension of an $R$-module $M$
is defined by $$\Rfd_RM:=\sup\{i\in\mathbb{N}_0|\Tor_i^R(M,T)\neq 0
\  \   \text{for some R-module} \ T \ \text{with finite flat
dimension} \}.$$

\begin{lemma} Let $\fa$ be an ideal of the Noetherian ring $R$ and
$Q$ a Gorenstein flat $R$-module. Then $\Rfd_R\Lambda^{\fa}(Q)=0$.
Moreover, if $R$ possesses a dualizing complex, then
$\Lambda^{\fa}(Q)$ is Gorenstein flat.
\end{lemma}

\begin{prf} Since $Q$ is Gorenstein flat, there exists an exact sequence
$$X=0\to Q\to F_{-1}\to F_{-2}\to \dots \to F_{-s}\to \cdots,$$ where each
$F_i$ is flat and $K_i:=\im(F_{-i}\to F_{-(i+1)})$ is Gorenstein
flat for all $i\geq 1$. Let $T$ be an $R$-module of finite flat
dimension $s$ 'say. By Corollary 2.3, the functor $\Lambda^{\fa}(-)$
is exact on the full subcategory of Gorenstein flat $R$-modules.
Hence, we have the following short exact sequences
$$0\to \Lambda^{\fa}(Q)\to \Lambda^{\fa}(F_{-1})\to
\Lambda^{\fa}(K_1)\to 0$$ and $$0\to \Lambda^{\fa}(K_i)\to
\Lambda^{\fa}(F_{-(i+1)})\to\Lambda^{\fa}(K_{i+1}) \to 0$$ for all
$i\geq 1$. By \cite[1.4.7]{B}, $\Lambda^{\fa}(F)$ is flat for all
flat $R$-modules $F$. Therefore, using the above short exact
sequences successively yields that
$$\Tor^{R}_j(\Lambda^{\fa}(Q),T)\cong
\Tor^{R}_{j+1}(\Lambda^{\fa}(K_1),T)\cong \dots \cong
\Tor^{R}_{j+s}(\Lambda^{\fa}(K_s),T)=0,$$ for all $j\geq 1$. This
shows that $\Rfd_R\Lambda^{\fa}(Q)=0$.

Now, assume that $R$ possesses a dualizing complex. By Theorem 2.5,
${\bf L}\Lambda^{\fa}(Q)\simeq \Lambda^{\fa}(Q)$, and so by
\cite[Theorem 5.10 b)]{CFH}, $\Gfd_{R}\Lambda^{\fa}(Q)$ is finite.
Thus, by \cite[Theorem 3.19]{H}, it turns out that
$\Gfd_R\Lambda^{\fa}(Q)=\Rfd_R\Lambda^{\fa}(Q)=0$, and so
$\Lambda^{\fa}(Q)$ is Gorenstein flat.
\end{prf}

Part i) of the following corollary strengthens \cite[Theorem 5.10
b)]{CFH}.

\begin{corollary} Let $R$ be a Noetherian ring possessing a
dualizing complex and $\fa$ an ideal of $R$.
\begin{enumerate}
\item[i)] Let $X\in \mathcal{D}_{\sqsupset}(R)$. Then $\Gfd_R{\bf L}
\Lambda^{\fa}(X)\leq \Gfd_RX$.
\item[ii)] Let $Y\in \mathcal{D}_{\Box}^f(R)$ be a  complex such
that either its projective or injective dimension is finite and
$X\in \mathcal{D}_{\Box}(R)$. Then $\sup {\bf
L}\Lambda^{\fa}(X\otimes_R^{{\bf L}}Y)\leq \Gfd_RX+\sup Y.$
\end{enumerate}
\end{corollary}

\begin{prf} i) It follows by Theorem 2.5 and Lemma 2.7.

ii) The inequality certainly holds if either $X$ has infinite
Gorenstein flat dimension or $Y$ is exact. So, assume that the
Gorenstein flat dimension of $X$ is finite and $Y$ is non-exact. By
i), the Gorenstein flat dimension of ${\bf L}\Lambda^{\fa}(X)$ is
finite, hence Lemma 2.6 and \cite[Theorem 3.5]{CFH} yield that
$$\begin{array}{ll}\sup {\bf L}\Lambda^{\fa}(X\otimes_R^{{\bf L}}Y)&=
\sup ({\bf L}\Lambda^{\fa}(X)\otimes_R^{{\bf L}}Y) \\
&\leq \Gfd_R{\bf L}\Lambda^{\fa}(X)+\sup Y\\
&\leq \Gfd_RX+\sup Y.
\end{array}
$$
\end{prf}

\begin{lemma} Let $R$ be a Noetherian ring possessing a dualizing
complex and $\fa$ an ideal of $R$.  The following are equivalent:
\begin{enumerate}
\item[i)]$\Lambda^{\fa}(Q)$ is flat for all Gorenstein flat $R$-module $Q$.
\item[ii)] $\Gfd_RQ=\fd_RQ$ for all $\fa$-adic complete $R$-modules
$Q$.
\end{enumerate}
\end{lemma}

\begin{prf} $i)\Rightarrow ii)$ Let $N$ be a $\fa$-adic complete $R$-module.
We have to show that $\Gfd_RN=\fd_RN$. Since, by \cite[Theorem
3.19]{H}, $\Gfd_RN\leq \fd_RN$, we can assume that $n:=\Gfd_RN$ is
finite. By \cite[Proposition 2.5]{Si} and its proof, one can choose
a flat resolution  $$F=\cdots \to F_i\overset{d_i}\to F_{i-1}\to
\cdots \to F_0\to 0$$ of $N$ such that $F_i$ and $\ker d_{i+1}$ are
$\fa$-adic complete for all $i\geq 0$. Let $Q:=\ker d_{n-1}$. Then
$Q$ is $\fa$-adic complete and by \cite[Theorem 3.14]{H}, $Q$ is
Gorenstein flat. Thus i) implies that $Q$ is flat, and so $N$ has
finite flat dimension. Now, use \cite[Theorem 3.19]{H} again to
deduce that $\Gfd_RN=\fd_RN$.

$ii)\Rightarrow i)$  Let $Q$ be a Gorenstein flat $R$-module. Then,
by Lemma 2.7, the $\fa$-adic complete $R$-module $\Lambda^{\fa}(Q)$
is Gorenstein flat. Hence, by ii), it turns out that
$\fd_R\Lambda^{\fa}(Q)=\Gfd_R\Lambda^{\fa}(Q)=0$, and so
$\Lambda^{\fa}(Q)$ is flat.
\end{prf}

Next, we present a characterization of regularity of Gorenstein
local rings.

\begin{corollary}  Let $(R,\fm,k)$ be a local Gorenstein ring. The
following are equivalent:
\begin{enumerate}
\item[i)] $\Lambda^{\fm}(Q)$ is flat for all Gorenstein flat
$R$-modules $Q$.
\item[ii)] $\Gfd_RQ=\fd_RQ$ for all $\fm$-adic complete $R$-modules
$Q$.
\item[iii)] $R$ is regular.
\end{enumerate}
\end{corollary}

\begin{prf}  i) and ii) are equivalent by  Lemma 2.9.

Since $R$  is Gorenstein, we see that $R$ is a dualizing complex of
$R$ and $k$ has finite Gorenstein flat dimension. Recall that $R$ is
regular if and only if the flat dimension of the $\fm$-adic complete
$R$-module $k$ is finite and if and only if the flat dimension of
any $R$-module is finite. Thus \cite[Theorem 3.19]{H} implies that
ii) and iii) are equivalent.
\end{prf}

We thank anonymous referee for suggesting the following example.

\begin{example} In the above corollary, the assumption of Gorensteiness of $R$
can not be deleted. To this end, let $(R,\fm)$ be a non-Gorenstein
local ring such that $\fm^2=0$, for examples of such rings see
\cite[Example 4.2]{CH}. Then by \cite[Proposition 4.3]{CH} any
Gorenstein flat  $R$-module is free, and so by \cite[1.4.7]{B},
$\Lambda^{\fm}(Q)$ is flat for all Gorenstein flat $R$-modules $Q$.
\end{example}


\end{document}